\documentclass[11pt]{article}
\usepackage{amsmath,amssymb,amsthm}
\theoremstyle{plain}
\newtheorem{Th}{Theorem}
\newtheorem{Cor}{Corollary}
\newtheorem{Prop}{Proposition}
\newtheorem{Lm}{Lemma}
\newtheorem{Def}{Definition}

\theoremstyle{definition}

\newtheorem{rem}{Remark}

\newcommand{\A}{{\cal A}}

\newcommand{\Z}{\mathbb{Z}}
\newcommand{\R}{\mathbb{R}}
\newcommand{\N}{\mathbb{N}}

\newcommand{\Q}{\mathbb{Q}}
\newcommand{\cc}{{\cal C}}
\newcommand{\ns}{\mbox{ns}}

\newcommand{\T}{{\cal T}}

\newcommand{\U}{{\cal U}}
\newcommand{\I}{{\cal I}}

\newcommand{\V}{{\cal V}}
\newcommand{\K}{{\cal K}}
\newcommand{\F}{{\cal F}}
\newcommand{\W}{{\cal W}}
\newcommand{\B}{{\cal B}}

\newcommand{\e}{\varepsilon}
\renewcommand{\l}{\lambda}
\newcommand{\f}{\varphi}

\renewcommand{\a}{\alpha}
\renewcommand{\b}{\beta}
\renewcommand{\d}{\delta}

\newcommand{\s}{\sigma}
\renewcommand{\c}{\gamma}
\newcommand{\sqe}{\subseteq}
\renewcommand{\l}{\lambda}

\newcommand{\la}{\langle}
\newcommand{\ra}{\rangle}
\renewcommand{\ns}{\mbox{ns}}
\newcommand{\inte}{\mbox{int}}
\newcommand{\range}{\mbox{range}}

\renewcommand{\*}{\,^*\!}
\renewcommand{\o}{\,^\circ\,\!\!}
\renewcommand{\vector}[3]{#1_{#2},\dots,\,#1_{#3}}
\newcommand{\all}{\forall}
\newcommand{\ex}{\exists}
\newcommand{\Liff}{\Longleftrightarrow}
\author{L.Yu.\ Glebsky, E.I.\ Gordon, C.\ Ward Henson\thanks{The research of
Gordon and Henson for this paper was partially supported by NSF
Grant DMS-9970009; Henson's research was also partially supported
by DMS-0100979 and DMS-0140677}}
\title{On finite approximations of topological algebraic systems}
\date{}
\usepackage{latexsym}
\begin{document}
\maketitle

\begin{abstract}
We introduce and discuss a concept of approximation of a
topological algebraic system $A$ by finite algebraic systems from
a given class $\K$.  If $A$ is discrete, this concept agrees with
the familiar notion of a \emph{local embedding} of $A$ in a class
$\K$ of algebraic systems. One characterization of this concept
states that $A$ is locally embedded in $\K$ iff it is a subsystem
of an ultraproduct of systems from $\K$. In this paper we obtain a
similar characterization of approximability of a locally compact
system $A$ by systems from $\K$ using the language of nonstandard
analysis.

In the signature of $A$ we introduce \emph{positive bounded}
formulas and their \emph{approximations}; these are similar to
those introduced by Henson \cite{he} for Banach space structures
(see also \cite{hi,hm}). We prove that a positive bounded formula
$\f$ holds in $A$ if and only if all precise enough approximations
of $\f$ hold in all precise enough approximations of $A$.

We also prove that a locally compact field cannot be approximated
arbitrarily closely by finite (associative) rings (even if the
rings are allowed to be non-commutative). Finite approximations of
the field $\R$ can be considered as possible computer systems for
real arithmetic. Thus, our results show that there do not exist
arbitrarily accurate computer arithmetics for the reals that are
associative rings.
\end{abstract}

\section{Introduction}

The numerical systems implemented in computers for simulation of
the field $\R$ are based on representation of reals in
floating-point form.  These systems are finite algebras with two
binary operations $\oplus$ and $\otimes$. The underlying set of
any such system $R$ is a finite, symmetric subset of $\R$ ($a\in
R$ iff $-a\in R$ for all $a\in\R$) on which the operations
$\oplus$ and $\otimes$ are defined as follows. Let $N$ be the
maximum of $R$. If $x,y\in R$ and $x+y\ \ (\mbox{resp., } x\times
y)\in [-N,N]$ then $x\oplus y\ \ (\mbox{resp., } x\otimes y)$ is
the element of $R$ nearest to $x+y\ \ (\mbox{resp., } x\times y)$.
Here $+$ and $\times$ are the addition and the multiplication in
$\R$. If $x+y\ \ (\mbox{resp., } x\times y)\notin [-N,N]$ then
$x\oplus y\ \ (\mbox{resp., } x\otimes y)$ is defined more or less
arbitrarily. If such overflow happens during a computation, the
numerical result might be incorrect; hence it is necessary to take
care that the overflow not occur. (In a floating-point system,
this is called \emph{exponent overflow}; see \cite[section
4.2.1]{kn}.)

The elements of a floating-point system $R$ are distributed
unevenly in the interval $[-N,N]$; they become especially sparse
when one gets close to the endpoints of this interval. This
non-uniformity entails a significant loss of accuracy in
calculations with large numbers, even if the results of
intermediate operations stay within the interval $[-N,N]$.
However, there exists an interval $[-M,M]\subseteq [-N,N]$ and a
positive $\e$ such that $R$ is $\e$-dense in $[-M,M]$ and for
every $x,y\in R$ if $x+y$ (resp., $x\times y$) $\in [-M,M]$, then
$x\oplus y$, (resp., $x\otimes y$) approximates $x+y$ (resp.,
$x\times y$) with an error that does not exceed $\e$.  By choosing
the parameters of the floating-point system correctly, $M$ can be
made arbitrarily large and $\e$ arbitrarily small.

According to the main definition in this paper (Definition 1),
this means that floating-point systems provide a family of
arbitrarily close finite approximations of the field $\R$
considered as a topological algebra. The algebraic properties of
systems based on the floating point representation are discussed
in [17], where it is shown that they are neither associative nor
distributive.

More generally, we consider in this paper ``continuous''
expansions of the field of real numbers; these are universal
algebras of the form $\mathcal{R} = (\R,1,\times,+,f_1,\dots,f_m)$
where the operations $f_j$ are continuous.  Several interesting
questions about the general nature of approximations of such
structures arise naturally.

First, is there a general procedure for constructing approximate
versions of theorems about continuous expansions of the reals? A
strong version of this question is the following: given a
proposition $\f$ about such structures, can one construct propositions
$\f_{M,\e}$ such that $\f$ holds for a given continuous expansion
$\mathcal{R}$ of the reals if and only if for all large enough $M$
and small enough $\e$, the proposition $\f_{M,\e}$ holds for all
finite systems $R$ approximating $\mathcal{R}$ on the interval
$[-M,M]$ with accuracy bounded by $\e$?  In section 4 we do
exactly this in an explicit way for positive first order sentences
$\f$ in which each quantifier is restricted to a bounded interval
of reals (Corollary \ref{pos_bound_sent} to Theorem
\ref{pos_bound}). It seems very difficult to do this in a more
general setting.

This kind of question may be important for an understanding of the
following type of problem. Suppose we use some convergent
numerical method for computation of a real function, or a
functional, or an operator. The theorem about convergence of this
method is a theorem about the field $\R$ but in our computer-based
``applications'' of this theorem we use a finite system $R$, which
only approximates $\R$. Can we be sure that the result of our
computation is approximately correct if we can use large enough
numbers and high enough accuracy? The fact that this problem is
natural can be demonstrated by the following example (concerning
the approximation of $\sin x$), which is discussed in
\cite[section 3.8]{md}.

Although the Taylor series for $\sin x$ converges for all $x$, the
approximate computation of $\sin x$  for large $x$ based on its
Taylor expansion gives an incorrect answer in a floating-point
system. For large $x$, the first few terms in a partial sum of
this series are also very large. Due to the fixed number of digits
in the floating-point representation of real numbers, the addition
of terms in a partial sum of the series should be done with the
terms taken in ascending order, to avoid roundoff error; this is
explained in \cite[chapter 2]{md}. However, calculation of the
$k^{th}$ term of the Taylor series for $\sin x$ produces exponent
overflow for large $x$ and $k$.

A second natural question concerning finite algebraic systems
approximating $\R$ is the following. What properties of continuous
expansions of the reals can hold for some finite systems that
approximate them arbitrarily closely? For example, let $\f$ be any
first order theorem about the field $\R$; is it true that for any
big enough $M$ and small enough $\e$ there exists a finite system
$R$ approximating $\R$ on the interval $[-M,M]$ with accuracy $\e$
such that $\f$ itself holds for $R$? We mentioned above that the
operations $\oplus$ and $\otimes$ in numerical systems based on the
floating-point representation are neither associative nor
distributive. Is it possible to construct finite \emph{rings} that
approximate $\R$ arbitrarily closely? (Here we answer this question
in the negative; see Theorem \ref{nonappr}.  It is easy to construct
approximating systems for $\R$ that are abelian groups for $\oplus$;
see Example 2 in section 2.)

These problems are discussed in the present paper in a more
general setting. We consider a locally compact algebraic system
$\A=\langle A,\theta\rangle$ of finite signature $\theta$ with
only function symbols (a universal algebra) and give a definition
of approximation of this system by a finite system $\A_f$ on a
compact set $C\subset A$ with accuracy $W$. Here $W$ is an element
of the uniformity on $A$ that defines its topology. We call $\A_f$
a $(C,W)$-approximation of $\A$. For example, if the topology on
$A$ is defined by a metric $\rho$, then we may take $W=\{\langle
x,y\rangle \in A^2\ |\ \rho(x,y)<\e\}$ for some $\e>0$. The
universal algebra $\A$ is said to be approximable by finite
algebras from a class $\K$ if for any $C$ and $W$ there exists a
$(C,W)$-approximation $\A_f\in\K$. The definition of approximation
of a locally compact group by finite groups discussed in \cite{go}
is a particular case of this definition. It is known \cite{go}
that all locally compact abelian groups are approximable by finite
groups but this is false in general for nonabelian groups
\cite{gr}. There exist groups that are approximable neither by
finite groups, nor by finite semigroups, nor even by finite
quasigroups \cite{agg} \cite{gr} \cite{gg}. It is proved in
\cite{gr} that the field $\R$ is not approximable by finite
fields; the signature here includes not only the operations of
addition and multiplication but also an operation giving the
multiplicative inverse of each nonzero element. Based on these
results we show here that locally compact fields are not
approximable by finite (associative) rings (Theorem 1). That is,
it is impossible to implement in a computer a numerical system for
arbitrarily accurate simulation of the field of reals that is a
finite (associative) ring.

In \cite{agk} (see also \cite{dhv}) finite approximations of
locally compact abelian groups are used for a construction of
finite dimensional approximations of pseudodifferential operators.
In this approach one simultaneously approximates the operators and
the group structures associated to them. This allows constructing
approximations which have nice properties (e.g., uniform
convergence and spectrum convergence). Usually, algebraic and
geometric structures connected with operators can be considered as
finite dimensional manifolds (e.g., the symmetry groups of
operators are often  Lie groups). Thus approximations of these
structures can be based on approximations of the field $\R$
together with some other continuous functions on $\R$.
Approximations of the other locally compact fields can be used in
$p$-adic analysis, adelic analysis, etc. This is another reason
for investigation of finite approximations of topological
algebraic systems.

Nonstandard analysis provides a natural language in which to
discuss approximate versions of statements about the reals; here
we return to the first general problem discussed in this
Introduction. For background on nonstandard analysis see, e.g.,
the recent books \cite{gol}, \cite{gkk}, and \cite{low}. A brief
introduction adequate for understanding sections 3 and 4 of this
paper is contained in \cite[Section 4.4]{kc}.

It is easy to construct approximate versions of first order
statements about continuous expansions $\mathcal{R}$ of the field
$\R$ using the language of nonstandard analysis, as we describe
next. Let $\f$ be a first order sentence in the language of
$\mathcal{R}$. Prenex rules and the presence of the arithmetic
operations $\times,+$ allow us to put $\f$ into an equivalent (in
$\mathcal{R}$) normal form
$$ Q_1x_1 \dots Q_mx_m\, [s = t] $$
where each $Q_j$ is either $\forall$ or $\exists$ and $s,t$ are
terms.  Now let $R$ be any hyperfinite approximation of
$\mathcal{R}$ (see Definition \ref{nonst-appr}) whose underlying
set is contained in $\* \R$; the mapping $j \colon R \to \* \R$ is
taken to be the inclusion.  It is then clear that $\f$ holds in
$\mathcal{R}$ if and only if the sentence
$$ Q_1^{fin}x_1 \dots Q_m^{fin}x_m\, [s \approx t] $$
holds in $R$; in a quantifier of the form $Q^{fin}x$ we take $x$
to range over the finite elements of $R$.  (See Proposition
\ref{nonst_form}.)

Standard reformulations of such nonstandard approximations can be
obtained using Nelson's algorithm \cite[Section 2]{ne} \cite{ne2}
for the translation of nonstandard statements into standard
language. Unfortunately, in the general case these standard
versions are extremely complicated. (Without using Nelson's
algorithm, we construct (section 4) comprehensible translations
for a large class of first-order sentences, the so-called positive
bounded sentences that we introduce here.)


Approximate versions of first-order sentences are discussed in
this paper for the general case of a locally compact algebra of
finite signature. The results obtained for our positive bounded
sentences are similar to well-known results about such sentences
in the theory of Banach spaces \cite{he} \cite{hh} \cite{hm} (see
also \cite{hi}).

The problem of constructing (nonstandard or standard) approximate
versions of higher order statements about $\R$ is also open and it
seems interesting and important. Solving it might lead to a deeper
understanding of the interaction between continuous mathematics
and its finite computer approximations.

The authors are grateful to the referee for valuable remarks and
important suggestions.

\section{Approximation of locally compact algebras}

Let $\A=\langle A,\theta\rangle$ be an algebraic system of finite
signature $\theta$ that contains only function symbols. We assume
that $A$ is endowed with a locally compact Hausdorff topology and
that the function symbols of $\theta$ are interpreted by
continuous functions.  (We denote these interpretations using the
same letters as the corresponding function symbols in $\theta$.)

Let $C\subset A$ be a compact set, $\U$ a finite covering of $C$
by relatively compact open sets (an r.c.o. covering),
$\A_f=\langle A_f,\theta\rangle$ a finite algebra of signature
$\theta$ and $j \colon A_f\to A$ a mapping. The interpretation of
a function symbol $g\in\theta$ in $A_f$ is denoted by $g_f$. For
$a_1,\dots a_n\in A_f$ we denote by $j(\langle
a_1,\dots,a_n\rangle)$ the $n$-tuple $\langle
j(a_1),\dots,j(a_n)\rangle$ . We say that $a,b\in C$ are
$\U$-close if $\exists U\in\U\ (a\in U\land b\in U)$.

\begin{Def} \label{c-u-appr}
\begin{enumerate}

\item We say that a set $M\subset A$ is a $(C,\U)$-grid
(equivalently, $M$ is a $\U$-grid for $C$) if for any $c\in C$ there
exists an $m\in M$ such that $c$ and $m$ are $\U$-close.

\item We say that $j$ is a $(C,\U)$-homomorphism if for any
$n$-ary function symbol $g\in\theta$ and for any $\bar a\in A_f^n$
such that $j(\bar a)\in C^n$ and $g(j(\bar a))\in C$, the elements
$g(j(\bar a))$ and $j(g_f(\bar a))$ are $\U$-close.

\item We say that the pair $\langle \A_f, j\rangle$ is a
$(C,\U)$-approximation of $\A$ if $j$ is a $(C,\U)$-homomorphism
and $j(A_f)$ is a $(C,\U)$-grid.

\item Let $\K$ be a class of finite algebras of signature
$\theta$. We say that the locally compact algebra $\A$ is
approximable by finite $\K$-algebras if for any compact set
$C\subset A$ and for any finite r.c.o. covering $\U$ of $C$ there
exists a $(C,\U)$-approximation $\langle \A_f,j\rangle$ of $\A$
such that $\A_f\in\K$.
\end{enumerate}
\end{Def}

\begin{rem} \label{discrete-appr}
If the topology on $A$ is discrete, then condition (4) in
Definition \ref{c-u-appr} is equivalent to the well-known
model-theoretic
 concept of local embedding of an algebraic
system $\la A,\theta\ra$ in a class $\K$ of algebraic systems of
the same signature $\theta$ (see e.g., \cite{ma}). The class of
discrete groups approximable by finite groups was studied in
\cite{vg}. It was shown, in particular, that in this case we
obtain the same class if we assume that the mapping $j$ is
injective. It is not known whether this is true for approximation
of topological algebras or even for approximation of discrete
algebras other than groups.
\end{rem}

\begin{rem} Note that if in the item 2 of Definition
\ref{c-u-appr} one has $\range(g)\cap C=\emptyset$ for all
$g\in\theta$ or $\range(j)\cap C=\emptyset$, then the mapping $j$
is a $(C,\U)$-homomorphism.
\end{rem}

Usually we deal with the case of a \emph{uniformly} locally compact
topology on $A$. This means that the topology on $A$ is determined
by a uniformity $\W$ and there exists  $W\in\W$ such that for any
$x\in A$ the set $W(x)=\{y\in A\ |\ \langle x,y\rangle \,\in W\}$ is
relatively compact. For example, all locally compact groups satisfy
this condition. For uniformly locally compact algebras of signature
$\theta$, we assume that the interpretations of function symbols are
continuous, but not necessary uniformly continuous. For example,
$\R$ is a uniformly locally compact space, but multiplication in
$\R$ is not uniformly continuous. It follows from the general theory
of uniform spaces (see, for example, \cite{bu}) that the restriction
of a continuous function to a compact subset $C$ is uniformly
continuous on $C$. For the case of uniformly locally compact
algebras Definition \ref{discrete-appr}(4) can be simplified.

We assume now that $\A$ is a uniformly locally compact algebra of
signature $\theta$ and $W$ is an element of the uniformity $\W$
such that $\forall x\in A\ [\overline{W(x)} \mbox{ is compact}]$.
(Here and below the closure of a set $E$ is denoted by $\overline
E$). Without loss of generality we may assume that $W$ is
symmetric (i.e., $\langle x,y\rangle \in W$ iff $\langle
y,x\rangle \in W$). The objects $C$, $\A_f$ and $j$ satisfy the
same assumptions as above. We say that $a,b\in C$ are $W$-close if
$\langle a,b\rangle \in W$.

\begin{Def} \label{c-w-appr}
\begin{enumerate}
\item We say that a set $M\subset A$ is a $(C,W)$-grid (equivalently,
$M$ is a $W$-grid for $C$) if for any $c\in C$ there exists an $m\in M$
such that $c$ and $m$ are $W$-close.

\item We say that $j$ is a $(C,W)$-homomorphism if for any $n$-ary
function symbol $g\in\theta$ and for any $\bar a\in A_f^n$ such
that $j(\bar a)\in C^n$ and $g(j(\bar a))\in C$, the elements
$g(j(\bar a))$ and $j(g_f(\bar a))$ are $W$-close.

\item We say that a pair $\langle \A_f, j\rangle$ is a
$(C,W)$-approximation of $\A$ if $j$ is a $(C,W)$-homomorphism and
$j(A_f)$ is a $(C,W)$-grid. If $A_f\subset A$ and $j \colon A_f\to
A$ is the inclusion map, we say that $\A_f$ is a
$(C,W)$-approximation of $\A$.

\item Let $\K$ be a class of finite algebras of signature
$\theta$. We say that a uniformly locally compact algebra $\A$ of
signature $\theta$ is approximable by finite $\K$-algebras if for
any compact set $C\subset A$ and for any $W\in\W$ such that
$\forall x\in A\ \overline {W(x)} \mbox{ is compact}$, there
exists a $(C,W)$-approximation $\langle \A_f,j\rangle$ of $\A$
such that $\A_f\in\K$.
\end{enumerate}
\end{Def}

We omit the simple proofs of the following four propositions.

\begin{Prop} \label{every}
For every uniformly locally compact algebra $\A$, any compact set
$C \sqe A$ and any element of the uniformity $W\in\W$ such that
$\forall x\in A\ [\overline{W(x)} \mbox{ is compact}]$, there
exists a finite $(C,W)$-approximation of $\A$.
\end{Prop}

\begin{Prop} \label{w-appr}
If $A$ is a compact set, then $\A$ is approximable by finite
$\K$-algebras in the sense of Definition \ref{c-w-appr} iff for
any $W\in\W$ there exists a finite $\K$-algebra that is an
$(A,W)$-approximation of $\A$.
\end{Prop}

\begin{Prop} \label {order1}
If $\la A_f,j\ra$ is a $(C,W)$-approximation of $\A$ in the sense
of Definition \ref{c-w-appr}, and we have a compact set
$C'\subseteq C$ and $\W\ni W'\supseteq W$, then the pair $\la
A_f,j\ra$ is a $(C',W')$-approximation of $\A$.
\end{Prop}

\begin{Prop} \label{equiv} A uniformly locally compact algebra $\A$ is
approximable by finite $\K$-algebras in the sense of Definition
\ref{c-u-appr} iff it is approximable by finite $\K$-algebras in
the sense of Definition \ref{c-w-appr}.
\end{Prop}

\begin{rem} Proposition \ref{equiv} shows that approximability of
a uniformly locally compact algebra $\A$ by finite $\K$-algebras is
a topological property; it does not depend on the uniformity on $\A$
but only on its topology. This fact is significant; for example, a
topology on a locally compact group is determined by the left
uniformity and by the right uniformity. It is well-known that these
two uniformities are not equivalent for some classical groups, e.g.,
for the group $SL(2,\R)$.

Since we deal only with uniformly locally compact algebras, in
what follows we use only Definition \ref{c-w-appr}, although all
of our results hold for the general case (after obvious
reformulations and modifications of proofs).
\end{rem}

Let $\A$ be a metric locally compact algebra with metric $\rho$ and
$W_{\e}=\{\la x,y\ra\ |\ \rho(x,y)<\e\}$.  In this situation we will
write $(C,\e)$-approximation (-grid, -homomorphism, etc) instead of
$(C,W_{\e})$-approximation (-grid, -homomorphism, etc).  Similarly,
we will write $\e$-grid for $C$ instead of $(C,W_{\e})$-grid.

Next we consider some examples of approximations of the field $\R$.
We use the signature $\s=\langle +,\times\rangle$. Since any compact
set $C\subset\R$ is contained in the interval $[-a,a]$ for some $a$,
and the sets $W_{\e}=\{\langle x,y\rangle \ |\ |x-y|<\e\},\ \e>0$
form a base of the uniformity on $\R$, it is enough to consider only
the $\left([-a,a], W_{\e}\right)$~-approximations of $\R$. We will
refer to them as $(a,\e)$-approximations. \label{a-epsilon}

\bigskip\noindent
{\bf Example 1} Recall that the floating-point form of a real
$\alpha$ is:
\begin{equation}
\alpha=\pm 10^p\times 0.a_1a_2\dots,
\end{equation}
where $p\in\Z$, and $a_1a_2\dots$ is a finite or infinite sequence
of digits such that $0\leq a_n\leq 9$, and $a_1\neq 0$. The
integer $p$ is called the exponent of $\alpha$, and
$0.a_1a_2\dots$, its normalized fraction part or mantissa.  Our
discussion of floating-point arithmetic mainly follows that of
\cite{kn}, differing only in some inessential technical details.

Fix natural numbers $P,Q$ and consider the finite set $A_{PQ}$ of
reals of the form (1) such that the exponent $p$ of $\alpha$
satisfies $|p|\leq P$ and the mantissa of $\alpha$ contains no
more than $Q$ decimal digits. We define binary operations $\oplus$
and $\otimes$ on $A_{PQ}$. In what follows, $ * $ stands for
either $+$ or $\times$. Let $\alpha,\beta\in A_{PQ}$ and suppose
the normal form of $\alpha * \beta$ is
\begin{equation}
\alpha * \beta=\pm 10^r\times 0.c_1c_2\dots\,.
\end{equation}
Note that the mantissa of $\alpha * \beta$ may contain more than
$Q$ digits.  In the following definition, the symbol $\circledast$
stands for $\oplus$ or for $\otimes$, depending on whether $*$
stands for $+$ or for $\times$. Then we define
\[
\alpha\circledast\beta=\left\{\begin{array}{ll} \pm 10^r\times
0.c_1c_2\dots c_Q\
&\mbox{if}\ |r|\leq P, \\
\pm 10^{P}\times0.\underbrace{99\dots9}_{Q\ \mbox{\scriptsize
digits}}\
&\mbox{if}\ r>P, \\
0\ &\mbox{if}\ r<-P.\end{array}\right.
\]
If the mantissa of $\alpha * \beta$ contains fewer than $Q$ digits
we complete it to a $Q$-digit mantissa by adding zeros at the
right.

Denote by $\A_{PQ}$ the algebra $\la A_{PQ},\oplus,\otimes\ra$ in
which the interpretations of the function symbols $+$ and $\times$
are the functions $\oplus$ and $\otimes$, respectively. It is easy
to see that for any positive $a$ and $\e$ there exist natural
numbers $P$ and $Q$ such that the algebra $\A_{PQ}$ is an
$(a,\e)$-approximation of $\R$. The systems $\A_{PQ}$ are
implemented in working computers. What properties of addition and
multiplication of the field of reals hold for $\oplus$ and
$\otimes$?

It is easy to see that the operations $\oplus$ and $\otimes$ are
commutative, $\xi\oplus\left(-\xi\right)=0$ and $\xi\oplus 0=\xi$
for any $\xi\in A_{PQ}$. Let $\alpha=\beta=0.60\dots 06$ and
$\gamma=0.60\dots 05$ (with $Q$ digits after the decimal point).
Then $\alpha\oplus\beta=\alpha\oplus\gamma$, so the cancellation
law fails for $\oplus$. Thus the associative law also fails for
$\oplus$. It is easy to construct examples to show that the laws
of associativity for $\otimes$ and distributivity for
$\oplus,\otimes$ in $\A_{PQ}$ also fail. See \cite[section
4.2.2]{kn} for some other identities of real arithmetic that hold
in these floating-point systems.

\bigskip\noindent
{\bf Example 2} Fix a natural number $M$ and a positive $\e$. Put
$A'_{M\e}=\{k\e\ |\ -M \leq k \leq M  \}$. Let $N=2M+1$. For any
$n\in\Z$ we will denote by $n$($\mbox{mod } N$) the unique element
of the set $\{ k \mid -M \leq k \leq M \}$ that is congruent to
$n$ modulo $N$. The operations $\oplus$ and $\otimes$ on
$A'_{M,\e}$ are defined as follows:
\begin{align}
k\e\oplus m\e&=(k+m)(\mbox{mod } N)\e, \\
k\e\otimes m\e&=[km\e](\mbox{mod } N)\e.
\end{align}

Denote by $\A'_{M,\e}$ the universal algebra in signature $\theta$
with the underlying set $A'_{M,\e}$ and the interpretation of the
function symbols defined by formulas (3) and (4). It is easy to
see that $\A'_{M\e}$ is an $(M\e,\e)$-approximation of $\R$. It is
obvious that $\A'_{M,\e}$ is an abelian group with respect to
$\oplus$ (see (3)). However, one can easily construct examples
which show that for any big enough $M$ and small enough $\e$ the
multiplication $\otimes$ satisfies neither the associative law nor
the distributive law.

\bigskip\noindent
{\bf Example 3} Consider approximation of the locally compact
field $\Q_p$ of $p$-adic numbers. Recall that any $p$-adic number
$\a\neq 0$ can be uniquely represented in the form
$$
\a=\sum\limits_{\nu=n}^{\infty}a_{\nu}p^{\nu},
$$
where $n \in \Z$ and for all $\nu \geq n$ in $\Z$ one has $0\leq
a_{\nu}<p$; moreover, the representation is normalized by taking
$a_n\neq 0$. The $p$-adic norm of $\a$ is then given by the
formula
$$
|\a|_p=p^{-n}, \eqno (8)
$$

The set $\Z_p=\{\a\ |\ |\a|_p\leq 1\}$ is a compact subring of
$\Q_p$, the ring of \emph{$p$-adic integers}. For any $m\in\Z$
consider the compact additive subgroup $p^{-m}\Z_p=\{\a\ |\
|\a|_p\leq p^{m}\}$. The sequence $\{p^{-m}\Z_p\ |\ m\in\Z\}$ is a
monotone sequence of compact sets that covers $\Q_p$. Hence, it is
enough to consider only the
$\left(p^{-m}\Z_p,p^{-n}\right)$-approximations of $\Q_p$ for all
$m,n \in \N$.

For any $n>0$, the set $p^n\Z_p$ is an ideal in $\Z_p$ and its
quotient ring $K_n$ is equal to $\Z/p^n\Z$. We represent an element
of this ring by its positive residue modulo $p^n$, so
$K_n=\{0,1,\dots, p^n-1\}$. We have $K_n\subset\Z_p$ as sets.
However, the ring operations in these sets are distinct. Indeed,
addition $+$ and multiplication $\times$ of natural numbers in
$\Z_p$ are the same as in $\N$, while addition $\oplus$ and
multiplication $\otimes$ in $K_n$ are equal to addition and
multiplication modulo $p^n$. For
$\a=\sum\limits_{\nu=0}^{\infty}a_{\nu}p^{\nu}\in\Z_p$ denote by
$\a_n$ the number $\sum\limits_{\nu=0}^{n-1}a_{\nu}p^{\nu}\in
{K_n}$. Then $|\a-\a_n|_p\leq p^{-n}$. Hence, $K_n$ is a
$p^{-n}$-grid for $\Z_p$. It is easy to see that
$$
|\a+\b-(\a_n\oplus\b_n)|_p,\ |\a\times\b-(\a_n\otimes\b_n)|_p\leq
p^{-n}.
$$

Thus, the inclusion map of $K_n$ into $\Z_p$ is a
$p^{-n}$-homomorphism. Hence, the ring $K_n$ is a
$(\Z_p,p^{-n})$-approximation of $\Z_p$. It follows that the
compact ring $\Z_p$ is approximable by finite commutative
associative rings. (See Proposition \ref{w-appr}.)

To construct a $(p^{-m}\Z_p,p^{-n})$-approximation ($0<m<n$) of
$\Q_p$ consider the set $H_{m,n}\subset p^{-m}\Z_p$ of all numbers
of the form $\sum\limits_{\nu=-m}^{n-1}a_{\nu}p^{\nu}$. Obviously,
$H_{m,n}$ is a $p^{-n}$-grid for $p^{-m}\Z_p$. We define operations
$\widehat\oplus$ and $\widehat\otimes$ such that the inclusion map
of $H_{m,n}$ into $p^{-m}\Z_p$ is a
$(p^{-m}\Z_p,p^{-n})$-homomorphism.

Note that $\a\in H_{m,n}$ iff $p^m\a\in K_{m+n}$. For any
$\a,\b\in H_{m,n}$ put
$$
\a\widehat\oplus\b=p^{-m}\left(p^m\a\oplus p^m\b\right),
$$
where $\oplus$ is the addition in $K_{m+n}$.

The definition of $\widehat\otimes$ is more complicated. Let
$$
p^m\a\times p^m\b=c_0+c_1p+\dots +
c_{m-1}p^{m-1}+c_mp^m\dots+c_{2m+2n-2}p^{2m+2n-2}.
$$
Put
$$
\a\widehat\otimes\b=c_mp^{-m}+\dots+c_{2m+n-1}p^{n-1}.
$$

It is easy to see that for all $\a,\b\in H_{m,n}$
$$|\a+\b-\a\widehat\oplus\b|_p\leq p^{-n},$$
and if $\a\times\b\in p^{-m}{\mathbb Z_p}$, then
$$|\a\times\b-\a\widehat\otimes\b|_p\leq p^{-n}.$$
Note that $\a\times\b\in p^{-m}{\mathbb Z_p}$ iff $c_k=0$ for
$k<m$.

This proves that the inclusion map of $H_{m,n}$ in $\Q_p$ is a
$(p^{-m}\Z_p,p^{-n})$-homomorphism.

Obviously, $\langle H_{m,n},\widehat\oplus\rangle$ is an abelian
group isomorphic to the additive group of $K_{m+n}$.

It is easy to see that for any integer $c$ such that $0\leq c<p$
one has $\frac c{p^m}\widehat\otimes\frac 1p=0$. Thus $\frac
1{p^m}\widehat\otimes\frac
1p\widehat\oplus\frac{p-1}{p^m}\widehat\otimes\frac 1p=0$, while
$\left(\frac
1{p^m}\widehat\oplus\frac{p-1}{p^m}\right)\widehat\otimes\frac
1p=\frac 1{p^m}.$ This shows that the distributive law fails for
$\widehat\oplus$ and $\widehat\otimes$.

Since $0\widehat\otimes p=0$ and $\frac 1p\widehat\otimes p=1$, we
have $\left(\frac 1{p^m}\widehat\otimes\frac
1p\right)\widehat\otimes p=0$, while $\frac
1{p^m}\widehat\otimes\left(\frac 1p\widehat\otimes p\right)=\frac
1{p^m}$. This shows that the associative law fails for
$\widehat\otimes$.

In all these examples, the finite algebras that approximate the
locally compact fields fail to be rings. Indeed, this is
inevitable, as the following theorem shows:

\begin{Th} \label{nonappr}
No infinite locally compact field can be approximated by finite
(associative) rings.
\end{Th}

\medskip\noindent
{\bf Proof}. Let $K$ be a locally compact field, $K_+$ the
additive group of $K$, and $K^{\times}$ the multiplicative group
of $K$. In this proof we denote the multiplication in $K$ by
$\cdot$. This multiplication is a continuous action of
$K^{\times}$ on $K_+$. It is obvious that this action does not
preserve the Haar measure on $K_+$.

Recall that a locally compact group is said to be unimodular if
the left and right Haar measures coincide.

It is well known \cite{lu} that if a unimodular group $G$ acts
continuously on a unimodular locally compact group $H$ by
automorphisms  and this action does not preserve the Haar measure
on $H$, then the semidirect product of $G$ and $H$ is
non-unimodular.

Thus, the semidirect product $K_+\leftthreetimes K^{\times}$ is a
non-unimodular group. This semidirect product is isomorphic to the
matrix group
$$
G=\left\{\left(\begin{array}{ll} a&b \\ 0&1\end{array}\right) \mid
\ a\in K^{\times},\ b\in K\right\}.
$$

Let us assume that $K$ is approximable by finite associative rings
and prove under this assumption that $G$ is approximable by finite
semigroups.

The group $G$ is homeomorphic to $K^{\times}\times K$ as a
topological space. Put
$$W_{\e}=\{\langle (a,b),(a',b')\rangle \ |\ a,a'\in K^{\times},\ b,b'\in K,
\max\{|a-a'|_K,\ |b-b'|_K\}<\e\},$$ where $|\cdot |_K$ is the norm
in $K$

We have to show that for any compact sets $A\subset K^{\times}$
and $B\subset K$ there exists a $(A\times B,
W_{\e})$-approximation $\langle S,j\rangle$ such that $S$ is a
semigroup.

Let $D=A\cup B\cup (A\cdot B)$. Then $D$ is a compact subset of
$K$ (here $A\cdot B=\{a\cdot b\ |a\in A,\ b\in B\}$). For any
positive $\d$ denote by $U_{\d}$ the set $\{\langle a,b\rangle \in
K\ |\ |a-b|_K<\d\}$ and put $C=\overline {U_{\e/2}(D)}$, where
$U_{\e/2}(D)=\bigcup\limits_{d\in D}U_{\e/2}(d)$. Since $D$ is a
compact set and any open ball in $K$ is relatively compact, we
have that $C$ is a compact set also.

According to our assumption, there exists a finite associative
ring $\langle F,\oplus,\odot\rangle$ and a map $j \colon F\to K$
such that the pair $\langle F,j\rangle$ is a
$(C,U_{\e/2})$-approximation of the field $K$. Our group $G$ is
equal to $K^{\times}\times K$ as a set. The multiplication in $G$
is given by the formula
$$
\langle a,b\rangle \cdot\langle c,d\rangle =\langle a\cdot
c,a\cdot d+b\rangle.
$$
Consider the finite set $S=F^*\times F$ and the multiplication in
$S$ given by the formula
$$
\langle s,p\rangle \odot \langle q,r\rangle =\langle s\odot
q,s\odot r\oplus p\rangle.
$$

Since $F$ is an associative ring, it is easy to see that $S$ is a
semigroup.

Define the map $i \colon S\to G$ by the formula
$$
i(\langle a,b\rangle)=\langle j(a),j(b)\rangle.
$$

Since $j(F)$ is an $U_{\e/2}$-grid for $C$ and thus, for $A$ and
$B$, it is obvious that $i(S)$ is a $W_{\e/2}$-grid for $A\times B$.

Let $i(\langle s,p\rangle),i(\langle q,r\rangle),i(\langle
s,p\rangle\odot\langle q,r\rangle)\in A\times B$. By the
definition of $i$, we have $j(s),j(q),j(s\odot q)\in A$. Hence,
$$
|j(s)\cdot j(q) -j(s\odot q)|_K<\frac{\e}2.
$$

Since $j(p),j(r)\in B$, we have $j(s)\cdot j(r)\in A\cdot B\subset
D\subset C$. Thus, $|j(s)\cdot j(r)-j(s\odot r)|_K<\frac{\e}2$ and
$j(s\odot r)\in C$. Hence,
$$
|j(s)\cdot j(r)+j(p)-\left(j(s\odot r)+j(p)\right)|_K<\frac{\e}2
$$
and
$$
|j(s\odot r)+j(p)-j\left((s\odot r)\oplus p\right)|_K<\frac{\e}2.
$$
Thus,
$$
|j(s)\cdot j(r)+j(p)-j\left((s\odot r)\oplus p\right)|_K<\e.
$$

This shows that $i$ is an $(A\times B, W_{\e})$-homomorphism.

Thus, the non-unimodular group $G$ is approximable by finite
semigroups.

By Theorem 4 of \cite{gg}, if a locally compact group is
approximable by finite semigroups, then it is approximable by
finite groups.

By Corollary 1 of Theorem 1 of \cite{gg}, if a locally compact
group $G$ is approximable by finite groups (indeed, even if only
by finite quasigroups), then $G$ is unimodular. This contradiction
completes the proof. \hfill $\Box$

\section{Nonstandard characterization of approximability}

We first recall some well-known notions and results from
nonstandard analysis. (See the books \cite{gol} \cite{gkk}
\cite{low} or the brief introduction in \cite[Section 4.4]{kc} for
necessary background.)  In this section, $\#(M)$ denotes the
cardinality of the set $M$.

Let $\mathbb{U}$ be a nonstandard universe and $\kappa$ an infinite
cardinal. Recall that $\mathbb{U}$ is $\kappa^+$-saturated if for
any family $\F$ of internal sets in $\mathbb{U}$ such that
$\#(\F)\leq\kappa$ and $\F$ satisfies the finite intersection
property one has $\bigcap\F\neq\emptyset$.

Let $\A=\la A,\theta\ra$ be as in the previous section, a
uniformly locally compact algebra of finite signature $\theta$. We
again assume that $\theta$ contains only function symbols and that
they are interpreted by continuous functions, which we denote by
the same letters as the respective function symbols.

Let $\l$ be the least infinite cardinal greater than the weight of
the topology on $\A$ and the weight of the uniformity on $\A$.
(The weight of a topology on $\A$ is  the minimal cardinality of a
base of this topology and the weight of a uniformity on $\A$ is
the minimal cardinality of a base of this uniformity).

\begin{Prop} \label{weight}
There exists a family $\cc_{\l}$ of compact subsets of $A$ with
the following properties:
\begin{itemize}

\item $\#(\cc_{\l})\leq\l$;

\item $\cc_{\l}$ is closed under finite unions;

\item  for any $C\in\cc_{\l}$ the interior $C^{\circ}$ of $C$ is
nonempty and $\bigcup\{C^{\circ}\ |\ C\in \cc_{\l}\}=A$.
\end{itemize}

\end{Prop}

\medskip\noindent
{\bf Proof} Let $\W_{\l}$ be a base of the uniformity on $A$ of
cardinality less or equal to $\l$. Without loss of generality we
assume that $\W_{\l}$ consists of elements $W$ such that for all
$x \in A$ the set $W(x)$ is open and relatively compact. Let $D$
be a dense subset of $A$ such that $\#(D)\leq\l$ (to obtain such
$D$ pick an element from each set in a base of the topology on $A$
of least cardinality). Take $\cc_{\l}$ to be the family of all
finite unions of the sets $\overline{W(d)}$, where $W\in\W_{\l}$
and $d\in D$. Then $\cc_{\l}$ satisfies the conditions of the
proposition. \hfill $\Box$

Throughout this section we deal with an arbitrary but fixed
$\l^+$-saturated nonstandard universe ${\mathbb U}$, with a fixed
family $\cc_{\l}$ satisfying Proposition \ref{weight} and with a
base $\W_{\l}$ of the uniformity $\W$ such that
$\#(\W_{\l})\leq\l$.

The nonstandard extension $\*\A$ of $\A$ is the algebraic system
$\la\*A,\theta\ra$, where any function symbol $f\in\theta$ is
interpreted in$\*\A$ by the nonstandard extension $\*f$ of the
operation $f$ in $A$. In what follows we omit the symbol $\*$ in
notations of $\*\A$-operations; i.e., we denote the
interpretations of a function symbol $f\in\theta$ in $\A$ and in
$\*\A$ by the same letter $f$.

For $\a,\b\in\*A$ we write $\a\approx\b$ if $\all W\in\W\ \la
\a,\b\ra\in\*W$. In this case we say that $\a$ and $\b$ are
infinitesimally close. Obviously, $\a\approx\b$ iff
$\la\a,\b\ra\in\*W$ holds for all $W\in\W_{\l}$. An element $\a$
of $\* A$ is called nearstandard if $\a\approx a$ holds for some
$a\in A$. Since $\W$ is a Hausdorff uniformity, if such an $a\in
A$ exists, then it is unique. In this case $a$ is called the
standard part or the shadow of $\a$ and it is denoted by $\o\a$.
If $\a$ is nearstandard, then $\b\approx\a$ iff for every open set
$U$ containing $\o\a$ one has $\b\in\*U$. Recall that for any
$a\in A$ the external set
$$
\mu(a)=\bigcap\{\*U\ |\ a \in U \mbox{ and } U\ \mbox{is
open}\}=\bigcap\limits_{W\in\W_{\l}}\*W(a).
$$
is called the monad of $a$.
\newcommand{\impl}{\longrightarrow}
\newcommand{\Impl}{\Longrightarrow}
\newcommand{\st}{\mbox{st}}

Denote the family of all compact subsets of $A$  by $\cc$.

For $B\sqe\*A$, we denote by $\ns(B)$ the set of all nearstandard
elements of $B$. It is well-known that a set $C\sqe A$ is compact
iff $\ns(\*C)=C$. Thus, since $A$ is a locally compact space, we
have
$$
\ns(\*A)=\bigcup\limits_{C\in\cc}\*C=\bigcup\limits_{C\in\cc_{\l}}\*C
\eqno (9)
$$

Let $f\in\theta$ and $f \colon A^n\to A$ for some standard natural
number $n$. Since $f$ is a continuous function, it is well-known
that for any $\bar\a,\bar\b\in (\ns(\*A))^n$ one has

$$
\bar\a\approx\bar\b\Impl f(\bar\a)\approx f(\bar\b) \eqno (10)
$$

Note that implication (10) holds for arbitrary
$\bar\a,\bar\b\in (\*A)^n$ iff the function $f$ is uniformly continuous.

The statements (9) and (10) and some of their obvious
modifications can be found in any of the books concerning
nonstandard analysis that were mentioned above.

Statement (10) implies the following:
\begin{Prop} \label{nearstand-subalg}\begin{enumerate}
\item The external set $\ns(\*A)$ is closed under
$\theta$-operations; i.e., $\la\ns(\*A),\theta\ra$ is a subalgebra
of $\*\A$. We denote this subalgebra by $\ns(\*\A)$.

\item The mapping $\st \colon \ns(\*\A)\to \A$ defined by the
formula $\st(\a)=\o\a$ is a surjective homomorphism of algebras
such that
$$
\st(\a)=\st(\b)\Liff \a\approx\b.
$$
Thus the equivalence relation $\approx$ restricted to $\ns(\*A)$
is a congruence relation on $\ns(\*\A)$ and the algebra
$\ns(\*\A)/\approx$ is isomorphic to $\A$.
\end{enumerate}
\end{Prop}

Let $\A_h$ be a hyperfinite algebra of signature $\theta$; i.e
$\A_h=\la A_h,\theta\ra$, where $A_h$ is a hyperfinite set and
every function symbol $f\in\theta$ is interpreted by an internal
function, which is denoted by $f_h$.

\begin{Def} \label{nonst-appr} Let $\A_h$ be a hyperfinite algebra
of signature $\theta$ and let $j \colon A_h\to \*A$ be an internal
mapping satisfying the following conditions:
\begin{enumerate}

\item $\all a\in A\,\ex b\in A_h\, j(b)\approx a$;

\item if $f\in\theta$ is an $n$-ary function symbol in $\theta$,
$\bar b\in A_h^n$ and $j(\bar b)\in (\ns(\*A))^n$, then
$j(f_h(\bar b))\approx f(j(\bar b))$.
\end{enumerate}
Then we say that the pair $\la\A_h,j\ra$ is a hyperfinite
approximation of the algebra $\A$.
\end{Def}

Assume that $\la\A_h,j\ra$ is a hyperfinite approximation of $\A$.

Put $(A_h)_b=j^{-1}(\ns(\*A))$. Elements of the set $(A_h)_b$ are
said to be \emph{feasible}.

For $b_1,b_2\in A_h$ put $b_1\sim b_2$ if $j(b_1)\approx j(b_2)$.

We call $\sim$ \emph{the indiscernibility relation}. If $b_1\sim
b_2$ we say that the elements $b_1$ and $b_2$ are
\emph{indiscernible}.

Definition \ref{nonst-appr} and Proposition \ref{nearstand-subalg}
imply immediately the following:

\begin{Prop} \label{hyperfin-subalg}
\begin{enumerate}

\item The external set $(A_h)_b$ is closed under
$\theta$-operations; i.e., $\la (A_h)_b,\theta\ra$ is a subalgebra
of $\A_h$ We denote this subalgebra by $(\A_h)_b$.

\item The mapping $\st\circ j=\iota \colon (\A_h)_b\to \A$ is a
surjective homomorphism of algebras such that for any $b_1,b_2\in
(A_h)_b$ one has
$$
\iota(b_1)=\iota(b_2)\Liff b_1\sim b_2.
$$
Thus the indiscernibility relation $\sim$ restricted to $(\A_h)_b$
is a congruence relation on $(\A_h)_b$ and the algebra
$(\A_h)_b/\sim$ is isomorphic to $\A$.
\end{enumerate}
\end{Prop}

Put $M=\cc_{\l}\times\W_{\l}$ and consider the partial order
$\leq$ on $M$ such that if $m_i=\la C_i,W_i\ra\in M,\ i=1,2$, then
$$
m_1\leq m_2\Liff C_1\supseteq C_2\wedge W_1\subseteq W_2.
$$

For $m\in M$ denote the set $\{m'\in M\ |\ m'\leq m\}$ by $M_m$.

Let $\la C,W\ra\in \*M$. By the transfer principle, the internal
set $C$ is $\*$-compact. We recall the meaning of this notion. Let
$\T$ be the topology on $A$. For any statement $P$, the
$\*$-version of $P$ is obtained by restricting all quantifiers to
internal sets. Any standard set $S$ involved in $P$ should be
replaced by its nonstandard extension $\*S$. Thus, an internal set
$C\sqe A$ is $\*$-compact if for any internal family $\U\sqe\*\T$
such that $C\sqe\bigcup\U$, there exists a hyperfinite subfamily
$\V\sqe\U$ such that $C\sqe\bigcup\V$.

Again let $\A_h$ be a hyperfinite algebra of signature $\theta$
and let $j \colon A_h\to\*A$ be an internal mapping. If $\la
C,W\ra\in\*M$, then we say that $\la\A_h,j\ra$ is a
$(C,W)-$approximation of $\*A$ if $\la\A_h,j\ra$ and $\*A$ satisfy
Definition \ref{c-w-appr}(3).

We say that a pair $\la C,W\ra\in\* M$ is infinitesimal if for any
$\la D,V\ra\in M$ one has $\la C,W\ra\leq\la\* D,\* V\ra$.
Obviously, if $\la C,W\ra\in\* M_{\l}$ is infinitesimal then
$C\supset\ns(\* A)$.

The following lemma is an immediate consequence of our assumption
that the nonstandard universe ${\mathbb U}$ is $\l^+$-saturated.

\begin{Lm} \label{about_M}
Let $N$ be an internal subset of $\* M$. The following statements
hold.
\begin{enumerate}
\item If for every $\la D,V\ra\in M$ one has $\la\* D,\* V\ra\in
N$, then there exists an infinitesimal element $\la C,W\ra\in N$;

\item if $N$ contains all infinitesimal elements of $\* M$, then
there exists an $m=\la D,V\ra\in M$ such that $\* M_m\sqe N$.
\end{enumerate}
\end{Lm}

\begin{Lm} \label{infin-appr}
A pair $\la\A_h,j\ra$ is a hyperfinite approximation of the
algebra $\A$  iff $\la\A_h,j\ra$ is a $(C,W)$-approximation of
$\*\A$ for some infinitesimal $\la C,W\ra\in\* M$.
\end{Lm}

\medskip\noindent
{\bf Proof} ($\Leftarrow$) Let $a\in A$. Since $j(A_h)$ is a $(
C,W)$-grid, there exists $b\in A_h$ such that $\la a,j(b)\ra\in
W$. Thus, for any $V\in\W$ one has $\la a, j(b)\ra\in\* V$; i.e.,
$j(b)\approx a$. Let $f\in\theta$ be an $n$-ary function symbol
and $f_h$ its interpretation in $\A_h$; take $\bar a\in
(A_h)_b^n$; i.e., $j(\bar a)\in\ns(\* A)^n\subset C^n$. Since $f$
is a continuous function, we have $f(j(\bar a))\in\ns(\* A)$.
Hence $f(j(\bar a))\in C$. By Definition \ref{c-w-appr} and the
transfer principle, $\la j(f_h(\bar a)),f(j(\bar a))\ra\in W$.
Thus, $j(f_h(\bar a))\approx f(j(\bar a))$. So $\la A_h,j\ra$ is a
hyperfinite approximation of $\A$.

($\Rightarrow$) Let $\la\A_h,j\ra$ be a hyperfinite approximation
of $\A$. Obviously, for any $\la D,V\ra\in M$ the pair
$\la\A_h,j\ra$ is a $(\* D,\* V)$-approximation of $\*\A$. By
Lemma \ref{about_M}(1), there exists an infinitesimal $\la C,W\ra$
such that $\la\A_h,j\ra$ is a $(C,W)$-approximation of $\* \A$.
\hfill $\Box$

\begin{Th} \label{nonst-charact}
A uniformly locally compact universal algebra $\A$ of finite
signature $\theta$ is approximable by finite algebras from a class
$\K$ iff there exist a hyperfinite algebra $\A_h=\la
A_h,\theta\ra\in\*\K$ and an internal mapping $j \colon A_h\to
\*A$ such that the pair $\la\A_h,j\ra$ is a hyperfinite
approximation of $\A$.
\end{Th}

\medskip\noindent
{\bf Proof} ($\Rightarrow$) Let $\A$ be approximable by finite
$\K$-algebras and let $\la C_0, W_0\ra$ be an infinitesimal
element of $\*M$. By the transfer principle, there exists a
hyperfinite algebra $\A_h\in\*\K$ and an internal mapping $j
\colon A_h\to A$ such that the pair $\la \A_h,j\ra$ is
$(C_0,W_0)$-approximation of $\*\A$. By Lemma \ref{infin-appr},
$\la \A_h,j\ra$ is a hyperfinite approximation of $\A$.

($\Leftarrow$) Let $\la \A_h,j\ra\in\*\K$ be a hyperfinite
approximation of $\A$. By Lemma \ref{infin-appr}, $\la
\A_h,j\ra\in\*\K$ is a $(C_0,W_0)$-approximation for some
infinitesimal $\la C_0,W_0\ra\in\*M$. Then by Proposition
\ref{order1} and the transfer principle, the pair $\la \A_h,j\ra$
is a $(\*C,\*W)$-approximation of $\*\A$ for all $\la C,W\ra\in
M$. By the transfer principle (used in the opposite direction),
for every $\la C,W\ra\in M$ there exists a finite
$(C,W)$-approximation of $\A$ that belongs to $\K$. \hfill $\Box$

\begin{Cor} \label{every-hyp}
For every uniformly locally compact algebra $\A$, there exists a
hyperfinite approximation of $\A$.
\end{Cor}

\medskip\noindent
{\bf Proof} Take $\K$ to be the class of all finite algebras of
signature $\theta$ and apply Theorem \ref{nonst-charact} and
Proposition \ref{every}.  \hfill $\Box$

\begin{rem} It follows from Definition \ref{c-u-appr} and
Proposition \ref{equiv} that Theorem \ref{nonst-charact} holds if
we take our nonstandard universe only to be $\nu^+$-saturated,
where $\nu$ is the weight of the topology on $A$.
\end{rem}

The topology on $(\A_h)_b/\sim$ induced by its isomorphism to $\A$
can be defined in terms of the triple $\la A_h, (A_h)_b, \sim\ra$.
We will now do this in a more general setting.

Recall that an external subset of a $\l^+$-saturated universe is
called a $\sigma$-set (respectively, a $\pi$-set) if it can be
represented by a union (respectively, an intersection) of a family
of internal sets of cardinality  $\leq\l$. Obviously $(A_h)_b$ is
a $\s$-subset of $A$, while $\sim$ is a $\pi$-set contained in
$A^2$.

The above considerations provide motivation for the following:

\begin{Def} \label{triple}
We say that a triple $\tau=\la T,T_b,\rho\ra$ is an abstract
nonstandard topological triple if $T$ is an internal set, $T_b\sqe
T$ is a $\s$-subset and $\rho$ is a $\pi$-equivalence relation on
$T$ such that for every $\a\in T_b$ the set $\rho(\a)=\{\b\in T\
|\ \la\a,\b\ra\in\rho \}$ is contained in $T_b$. We call $T_b$ the
set of abstractly feasible elements and $\rho$ the abstract
indiscernibility relation. If $T$ is hyperfinite, we call $\tau$ a
hyperfinite topological triple.
\end{Def}

We now introduce a topology on the quotient set $\widehat
T=T_b/\rho$. For $\a\in T$ denote by $\a^{\rho}$ the
$\rho$-equivalence class of $\a$.

Let $F\subset T$. Put $i(F)=\{\a\in F\ |\ \a^{\rho}\subset F\}$.
Denote by $\I$ the family of all internal subsets of $T_b$. Let
$\T_{\tau}$ be the topology on $\widehat T$ obtained by taking the
family $\{i(F)^{\rho}\ |\ \a\in i(F),\,F \in\I\}$ to be a base of
neighborhoods of the point $\a^{\rho}$, for each $\a\in T_b$. Here
$i(F)^{\rho}=\{\b^{\rho}\ |\ \b^{\rho}\subseteq F\}$.

The construction of the topological space $(T_b,\T_{\tau})$ is a
generalization of the well-known construction of the nonstandard
hull. This generalization was introduced in \cite{aam} for the
case of hyperfinite abelian groups (see also \cite{go}).

\begin{Th} \label{abstr_appr}
\begin{enumerate}
\item The weight of the topology $\T_{\tau}$ is $\leq\l$.

\item The topological space $(\widehat T,\T_{\tau})$ is locally
compact iff for every internal set $F\subseteq T_b$ and for every
internal set $G$ such that $\rho\subseteq G\subseteq T\times T$
there exists a set $K\subseteq F$ of standard finite cardinality
that satisfies the following condition:
$$
F\subseteq\bigcap\limits_{\a\in K}G(\a).
$$

\item If $\f \colon T^n\to T$ is an internal $n$-ary operation on
$T$ for some standard $n$, and we assume that the set $T_b$ of
feasible elements of $T$ is closed under $\f$ and
$\f\upharpoonright T_b$ is stable under the indiscernibility
relation $\rho$; i.e.,
$$
\all \bar a,\bar a'\in T_b^n\,(\bar a\rho\bar a'\impl \f(\bar
a)\rho\f(\bar\a')),    \eqno (11)
$$

then the induced $n$-ary operation $\f^{\#}$ on $\widehat T$
(i.e., $\f^{\#}$ is such that for every $\bar a\in T_b^n$ one has
$\f^{\#}(\bar a^{\rho})=\f(\bar a)^{\rho}$) is continuous in the
topology $\T_{\tau}$.

\item Let $\la A_h,j\ra$ be a hyperfinite approximation of a
uniformly locally compact algebra $\A$, let $(\A_h)_b$ and $\sim$
be as defined in Proposition \ref{hyperfin-subalg}, let $\tau=\la
A_h, (A_h)_b,\sim\ra$, and $\hat\A_h=(\A_h)_b/\sim$. Then the
isomorphism of algebras $\hat\A_h$ and $\A$ induced by the
homomorphism $\iota \colon (\A_h)_b\to\A$ of Proposition~
\ref{hyperfin-subalg} is an isomorphism of topological algebras
with respect to the topology $\T_{\tau}$ on $\hat\A_h$.
\end{enumerate}
\end{Th}

A proof of this theorem for the case of locally compact abelian
groups is contained in \cite{aam} and in \cite{go}. It can be
transferred without any changes to the general case.

Let $\A_h=\la A_h,\theta\ra$ be an internal algebra, let
$(\A_h)_b=\la (A_h)_b,\theta\ra$ and let $\rho$ be a
$\pi$-equivalence relation on $A_h$. We say that the triple
$\tau=\la\A_h,(\A_h)_b,\rho\ra$ is a nonstandard topological
$\theta$-triple, if $\rho\upharpoonright(A_h)_b$ is a congruence
relation on $(\A_h)_b$ (i.e., (11) holds for all operations $\f$
from $\theta$) and $\la A_h, (A_h)_b,\rho\ra$ is an abstract
nonstandard topological triple.

Theorem \ref{abstr_appr}(2) shows that if
$\tau=\la\A_h,(\A_h)_b,\rho\ra$ is a nonstandard topological
$\theta$-triple, then the quotient algebra
$\hat\A_h=(\A_h)_b/\rho$ is a topological algebra with respect to
the topology $\T_{\tau}$.

We say that a topological algebra $\A=\la A,\theta\ra$ is
\emph{abstractly approximable by finite algebras from a class
$\K$}, if there exists a hyperfinite topological $\theta$-triple
$\tau=\la\A_h,(\A_h)_b,\rho\ra$ such that $\A_h\in\*\K$ and $\A$
is topologically isomorphic to $\hat\A_h$.

Theorem \ref{abstr_appr} together with Proposition
\ref{hyperfin-subalg} show that if $\A$ is approximable by finite
$\K$-algebras, then it is abstractly approximable by finite
$\K$-algebras.

The following question is open. Is is true that any locally
compact algebra $\A$ that is abstractly approximable by finite
$\K$-algebras is approximable by finite $\K$-algebras in the sense
of Definition \ref{c-u-appr}?

It is easy to see that Theorems 1 and 4 of \cite{gg} stay true if
we replace approximability (of groups by finite quasigroups and
finite semigroups) by {\it abstract} approximability (of groups by
finite quasigroups and finite semigroups).

This implies the following proposition, which strengthens Theorem
\ref{nonappr}.

\begin{Prop} \label{nonst-non-appr} No infinite locally compact field is
abstractly approximable by finite (associative) rings.
\end{Prop}

An interesting discussion about the relations between real
analysis and discrete analysis is contained in \cite{Ze}. The main
idea of that paper is expressed as follows: ``Continuous analysis
and geometry are just degenerate approximations to the discrete
world \dots . While discrete analysis is conceptually simpler
\dots\ than continuous analysis, technically it is usually much
more difficult. Granted, real geometry and analysis were necessary
simplifications to enable humans to make progress in science and
mathematics \dots''.

The discussion in this section shows how the idea that continuous
mathematics is an approximation of the discrete could be
formalized. We may assume that we deal only with finite sets, but
some of these sets are so big that they contain some only vaguely
defined subclasses, which do not satisfy all the properties of
sets. For example, the induction principle fails for these
subclasses. For example, recall the well-known paradox of the pile
of sand, due to Eubulides, IV century B.C.: one grain of sand is
not a pile, and if $n$ grains of sand do not form a pile, then
$n+1$ grains also do not form a pile; so, how can we get a pile of
sand? According to our approach, hyperfinite sets of infinite
cardinality simulate such large sets and external subsets simulate
their vaguely defined subclasses. This follows from the obviously
true statement: ``A hyperfinite set has a standard cardinality iff
all its subsets are internal''. Under this approach the set of all
grains of sand is hyperfinite and a pile of sand is an external
subset of this set\footnote{P.Vopenka \cite{vop} suggested an
axiomatic set theory for the formalization of this idea. The main
defect of his approach is its opposition to classical mathematics.
Another axiomatization of hyperfinite sets was suggested in
\cite{ang}, where classical mathematics was interpreted in the
framework of hyperfinite sets.}.

According to Proposition \ref{hyperfin-subalg} and Corollary
\ref{every-hyp}, for every locally compact algebra $\A$ there
exists a hyperfinite algebra $\A_h$, an external subalgebra
$(\A_h)_b$ and an equivalence relation $\sim$ such that $\A$ is
isomorphic to $(\A_h)_b/\sim$. So $(A_h)_b$ can be viewed as a
subclass of feasible elements and $\sim$ as an indiscernibility
relation.

Proposition \ref{nonst-non-appr} together with the results on
non-approximability of Lie groups from \cite{agg} explain, in some
sense, why continuous analysis is simpler than discrete analysis.
The discrete algebraic structures that are used in science need
not have algebraic properties as good as those possessed by the
corresponding continuous structures.

We complete this section with a formulation of the concept of
approximability in terms of ultraproducts.

If an algebra $\A$ is discrete (see Remark \ref{discrete-appr}),
then Definition \ref{c-w-appr}(4) of the concept of approximation
of $\A$ by finite $\K$-algebras can be reformulated in the
following way.

\begin{Prop} \label{d-appr}
A discrete algebra $\A$ of a finite signature $\theta$ is
approximable by finite $\K$-algebras iff for any finite subset
$C\subset A$ there exists a finite algebra $\A_C\in\K$ and a map
$j \colon A_C\to A$ such that
\begin{enumerate}
\item $C\subset j(A_C)$.

\item For any $n$-ary function symbol $f\in\theta$ and for any
$\bar a\in A_C^n$ such that $j(\bar a)\in C$ and $f(j(\bar a))\in
C$, one has
$$
j(f_C(\bar a))=f(j(\bar a)),
$$
where $f_C$ is the interpretation of $f$ in $\A_C$.
\end{enumerate}
\end{Prop}

In \cite{vg} are presented some examples of locally compact groups
$G$ such that $G$ is approximable by finite groups as a discrete
group but $G$ is not approximable by finite groups as a
topological group.

The following proposition is contained in \cite{ma}. (A proof can
also be found in \cite{agg}.)

\begin{Prop} \label{dis-ult}
A discrete algebra $\A$ is approximable by finite $\K$-algebras
iff $\A$ is isomorphic to a subalgebra of an ultraproduct of
finite $\K$-algebras.
\end{Prop}

Theorem \ref{nonst-charact} together with Proposition
\ref{hyperfin-subalg} can be considered as a generalization of
Proposition \ref{dis-ult} to the setting of approximation of
topological algebras.

Indeed, if our nonstandard universe is a $\l^+$-saturated
ultrapower of a standard universe, then any hyperfinite algebra
$\A_h\in\*\K$ is isomorphic to an ultraproduct of finite
$\K$-algebras. Internal subsets of $A_h$ correspond to subsets of
this ultraproduct that are ultraproducts themselves. Unions
(respectively, intersections) of at most $\l$ many internal
subsets are called $\s$-sets (respectively, $\pi$-sets). Combining
Theorem \ref{nonst-charact} and Proposition \ref{hyperfin-subalg}
with these remarks we obtain the following:

\begin{Prop} \label{top-ult} If a uniformly locally compact algebra
$\A$ of signature $\theta$ is approximable by finite algebras from
a class $\K$, then $\A$ is isomorphic to a quotient algebra of a
$\s$-subalgebra $\B_{\s}$ of some $\l^+$-saturated ultraproduct
$\B$ of finite $\K$-algebras with respect to some
$\pi$-equivalence relation $\rho$ on $\B$, such that
$\rho\upharpoonright\B_{\s}$ is a congruence relation. (Here $\l$
is the weight of the topology on $\A$).
\end{Prop}

The necessity in Proposition \ref{dis-ult} is a special case of
Proposition \ref{top-ult}.  Indeed, it is easy to see that if the
topology on $\A$ is discrete, then the equivalence relation $\rho$
is the relation of equality. (See Proposition
\ref{hyperfin-subalg}.)

\bigskip

\section{Positive bounded formulas and finite approximations}

In this section we consider first order statements true of a
locally compact algebra $\A$, and we investigate approximate
versions of those statements that hold in finite approximations of
$\A$. We start with approximations of statements that are
formulated in the language of nonstandard analysis.

Let $L_{\theta}$ be the set of all first order formulas in the
signature $\theta$ and let $\f(x_1,\dots,x_n)\in L_{\theta}$.
Denote by $\f_{\sim}$ the formula obtained from $\f$ by the
replacement of each atomic subformula $t_1=t_2$ by the formula
$t_1\sim t_2$; here $t_1$ and $t_2$ are terms in the signature
$\theta$,

Let $\la\A_h,j\ra$ be a hyperfinite approximation of $\A$. Then
the formula $\f_{\sim}$ has an obvious interpretation in the
algebra $(\A_h)_b$ of feasible elements of $\A_h$ (cf. Proposition
\ref{hyperfin-subalg}). Every term $t(x_1,\dots, x_n)$  of
signature $\theta$ is interpreted by a function $t_b$ on
$(A_h)_b^n$, obtained by substitution of the function $f_b$ for
any function symbol $f$ involved in $t$ (we denote the restriction
of $f_h$ to $(A_h)_b$ by $f_b$). Then for any $a_1,\dots, a_n\in
(A_h)_b$ one has $(\A_h)_b\models t_1(a_1,\dots,a_n)\sim
t_2(a_1,\dots, a_n)$ iff the elements $(t_1)_b(a_1,\dots,a_n)$ and
$(t_2)_b(a_1,\dots,a_n)$ are indiscernible. The following
proposition is an immediate corollary of Proposition
\ref{hyperfin-subalg}.

\begin{Prop} \label{nonst_form}
If $\la A_h,j\ra$ is a hyperfinite approximation of $\A$ then for
any formula $\f(x_1,\dots,x_n)\in L_{\theta}$ and any
$a_1,\dots,a_n\in (A_h)_b$ one has
$$
(\A_h)_b\models\f_{\sim}(a_1,\dots,a_n)\Longleftrightarrow
\A\models\f(\o j(a_1),\dots,\o j(a_n)).
$$
\end{Prop}

\begin{rem} \label{abstract} The same proposition is true also
for any abstract hyperfinite approximation of $\A$ if we replace
$\A$ by $\widehat \A_h$ and $\o j(a_i)$ by the canonical image of
$a_i$ in $\widehat\A_h$.
\end{rem}

>From the point of view of computer numerical systems discussed in
the Introduction, Proposition 2 has the following interpretation.
In the setting of nonstandard analysis, we can consider an
idealized computer that has a hyperfinite memory. Then the
numerical system $\R_h$ for simulating the field of reals that is
implemented in this computer is a hyperfinite algebra in the
signature $\s=\la +,\times\ra$ and $\R_h$ is a hyperfinite
approximation of $\R$. So Proposition \ref{nonst_form} provides a
lot of information about $\R_h$.

Suppose $N = \max \{ |\alpha| \mid \alpha \in \R_h \}$.  Then the
elements of $(\R_h)_b$ can be considered as elements that are far
enough from the end points of the interval $[-N,N]$ so that
exponent overflow never occurs in computations involving them. It
is very natural that the property of being ``far enough from the
end points of the interval $[-N,N]$'' is an external property: if
a natural number $n$ is `far enough from the end points'' then
obviously the same is true for $n+1$. Thus the induction principle
fails for this property. Proposition \ref{nonst_form} shows that
the first order properties of $\R$ hold approximately for the
computer implementation of $\R$, as long as we only consider
elements that are far enough from the end points of the interval
$[-N,N]$. This fact seems to be very clear for those who use
computers for numerical computation. The language of nonstandard
analysis makes it possible to formulate a rigorous mathematical
theorem that expresses this phenomenon.

\bigskip\noindent
{\bf Example 4} Consider the algebra $\A_{PQ}$ discussed in
Example 1 of section 2. It is easy to see that if
$P,Q\in\*\N\setminus\N$, then $\A_{P,Q}$ is a hyperfinite
approximation of $\R$ (here $j$ is the inclusion map). Consider a
formula $\f(x,y)$ of the signature $\s=\la +,\times\ra$. Let
$\R\models\all x\ex y\f(x,y)$. Put $\psi(x)=\ex! y\f(x,y)$,
$\eta(x)= \ex y_1,y_2(y_1\neq y_2\wedge\f(x,y_1)\wedge\f(x,y_2))$,
where $\ex! y\f(x,y)$ means that there exists a unique $y$ such
that $\f(x,y)$. Assume that for every rational number $\a$ one has
$$
\R\models\psi(\a) \eqno (12)
$$
Thus, $\*\R\models\psi(a)$ holds for every $a \in A_{PQ}$. Let us
assume also that there exists an irrational $\a$ such that
$$
\R\models\eta(\a) \eqno (13)
$$

Consider the following question. Given an arbitrary $\a\in\R$, how
can we determine whether $\a$ satisfies (12) or (13) using only
our computer? The qualitative answer to this question is the
following. If $\a$ satisfies (12), then for all precise enough
approximations $a_1$ and $a_2$ of $\a$, any $b_1$ and $b_2$ such
that $\f(a_1,b_1)$ and $\f(a_2,b_2)$ are true with a high accuracy
must be very close to each other. If $\a$ satisfies (13), then
there exist two arbitrarily precise approximations $a_1$ and $a_2$
of $\a$ and two significantly distinct $b_1$ and $b_2$ such
$\f(a_1,b_1)$ and $f(a_2,b_2)$ are approximately true.

A rigorous mathematical statement that reflects this qualitative
answer follows from Proposition \ref{nonst_form}. Indeed, it is
easy to see that
$$
(\A_{PQ})_b\models\all
x_1,x_2(\psi_{\sim}(x_1)\wedge\psi_{\sim}(x_2)\wedge x_1\sim
x_2\to\all y_1,y_2(\f_{\sim}(x_1,y_1)\wedge\f_{\sim}(x_2,y_2)\to
y_1\sim y_2)) \eqno (14)
$$

and
$$
(\A_{PQ})_b\models\all x(\eta_{\sim}(x)\to\ex
x_1,x_2,y_1,y_2(\f_{\sim}(x_1,y_1)\wedge\f_{\sim}(x_2,y_2)\wedge
x_1\sim x\wedge x_2\sim x\wedge\neg(y_1,\sim y_2))) \eqno (15)
$$

Let us illustrate this discussion by a very simple numerical
example. Consider the following system
$$
\left\{\begin{array}{rll} x+ay &=&b \\ ax+b
y&=&2\end{array}\right. \eqno(16)
$$

This system has
\begin{enumerate}
\item a unique solution, if $a^2\neq b$,

\item no solutions if $a^2=b\neq\sqrt[3]4$,

\item infinitely many solutions if $a^2=b=\sqrt[3]4$.
\end{enumerate}
In the last case the general solution is given by the formula
$$
x+\sqrt[3]2\cdot y=\sqrt[3]4  \eqno(17)
$$

Performing numerical calculations on a computer, we deal only with
rational numbers. Thus, the third case cannot occur in computer
calculations.

Taking the 5-digit approximations to $\sqrt[3]2$ and $\sqrt[3]4$
for $a$ and $b$ and solving the system (16) on a computer we
obtain the solution $x' = 0.74552,\ y' = 0.6682$, which satisfies
(17) with accuracy $10^{-5}$

Taking the 10-digit approximations to $\sqrt[3]2$ and $\sqrt[3]4$
for $a$ and $b$ and solving the system (16) we obtain the solution
$x' = -0.9979450387,\ y' = 2.051990552$, which satisfies (17) with
accuracy $10^{-10}$. We see that these two approximate solutions
of the system (16) are significantly distinct (compare with (15)).

In the language of nonstandard analysis it is only possible to
formulate mathematical theorems that give us some qualitative
picture of the connection between continuous problems and their
computer simulations. To obtain specific estimates it is necessary
(but not sufficient) to formulate a standard version of
Proposition \ref{nonst_form}.

In the language of classical mathematics, we can only consider
approximate properties of reals that hold eventually when the
memory of computers increases to infinity and the accuracy becomes
more and more precise. We will see that only a restricted result
can be obtained in this way.

We say that a formula $\f\in L_{\theta}$ is {\it positive} if it
can be built up from atomic formulas using only conjunctions,
disjunctions and quantifiers. The main result of this section
concerns positive formulas in prenex form
$$
Q_1y_1\dots Q_my_m\psi(\vector x1n, \vector y1m), \eqno(18)
$$
where the $Q_i$ are quantifiers and $\psi$ is a disjunction of
conjunctions of atomic formulas.

An arbitrary (not necessary positive) formula $\f$ is equivalent
to a formula in the form (18), where $\psi$ is a disjunction of
conjunctions of atomic subformulas of $\psi$ and negations of
atomic subformulas of $\psi$. Let $\vector \c1k$ be the list of
all atomic formulas and their negations involved in $\psi$. For
any $1\leq i\leq k$ fix  $W_i\in\W$ and denote by $\c_i[W_i]$ the
formula $\la t_1,t_2\ra\in W_i$ if $\c_i$ is $t_1=t_2$ and the
formula $\la t_1,t_2\ra\notin W_i$ if $\c_i$ is $\neg(t_1=t_2)$.
Here $t_1$ and $t_2$ are terms in the signature $\theta$.

Define the interpretations of the formula $\la t_1,t_2\ra\in W$ in
$\A$ and in an arbitrary $(C',W')$-approximation $\la\A_f, j_f\ra$
of $\A$, where $\la C', W'\ra\in M$, as follows:

Let $\tau_1,\tau_2\in A$ be interpretations of the terms $t_1$ and
$t_2$ in $\A$. Then $\A\models \la t_1, t_2\ra\in W$ iff
$\la\tau_1,\tau_2\ra\in W$.

If  $\xi_1$ and $\xi_2$ are interpretations of the terms $t_1$ and
$t_2$ in $\A_f$, then $\A_f\models\la t_1,t_2\ra\in W$ iff $\la
j_f(\xi_1),j_f(\xi_2)\ra\in W$.

Denote by $\f[\vector W1k]$ the formula that is obtained from $\f$
by replacement of each $\c_i$ by $\c_i[W_i]$ respectively. The
formula $\f[\vector W1k]$ is called an approximation of $\f$.
Obviously, if $\f$ is positive, then for any $\vector W1k\in\W$
one has $\f\Longrightarrow\f[\vector W1k]$ (for both
interpretations). This is not true in general if $\f$ is
non-positive. Similarly, if $W'_i\sqe W_i,\ i=1,\dots k$, then
$\f[\vector {W'}1k]\Longrightarrow\f[\vector W1k]$ for a positive
$\f$ but not generally if $\f$ is non-positive. For positive
formulas we say in this case that $\f[\vector {W'}1k]$ is a finer
approximation than $\f[\vector W1k]$. Obviously, for any
approximation $\f[\vector W1k]$ of a positive formula $\f$ there
exists a finer approximation $\f[\vector {W'}1k]$ such that
$W'_1=\dots=W'_k=W$ (it is enough to put $W=W_1\cap\dots\cap
W_k$). In this case we write $\f[W]$ instead of $\f[W,\dots,W]$.
In what follows we deal only with approximations of the form
$\f[W]$ of a positive formula $\f$.

If $B\sqe A$ and $Q$ is either $\forall$ or $\exists$ then
$Q_Bx\dots$ is interpreted in $\A$ by $\forall x(x\in B\to\dots)$
or $\exists x(x\in B\land\dots)$ and in a finite
$(C,W)$-approximation $\langle \A_f,j_f\rangle$ of $\A$ by
$\forall x(x\in j_f^{-1}(B)\to\dots)$ or $\exists x(x\in
j_f^{-1}(B)\land\dots)$. \label{ibq}

Quantifiers of the form $Q_B$ are called bounded quantifiers. If
all quantifiers in a formula $\f$ are bounded then we say that
$\f$ is bounded.

Let $c=\langle \vector C1m\rangle$ be an $m$-tuple of subsets of
$A$ and let $\f$ be a positive prenex formula as in (18). Then $\f[c]$
is the formula
$$
Q_{1C_1}y_1\dots Q_{mC_m}y_m\psi.    \eqno (19)
$$

A formula of the form (19) is said to be a \emph{positive bounded
formula}.

In what follows we consider only positive bounded formulas $\f[c]$
that satisfy the following condition:

for any $i\leq m$ such that $Q_i = \forall$ (respectively,
$Q_i=\exists)$ the set $C_i$ is a relatively compact open
(respectively, compact) set.

In this case we say that an $m$-tuple $c$ of subsets of $A$ is
$\f$-regular.

\bigskip\noindent
{\bf Example 5} Consider the signature $\s'$ obtained from the
signature $\s=\la\oplus,\otimes\ra$ by adding a constant for each
real number. Let $\f$ be a formula of the form (18) in the
signature $\s'$ and $c$ be a $\f$-regular $m$-tuple that consists
only of open and closed intervals. Since the relation $x\leq y$ is
expressed by the positive formula $\ex z(y=x+z^2)$ and the
universal quantifiers are restricted to open intervals, while the
existential quantifiers are restricted to closed intervals, it is
easy to see that $\f[c]$ is equivalent to a positive formula of
the signature $\s'$.

For two $\f$-regular $m$-tuples $c$ and $c'$ we say that $c\ll c'$
if for any $i\leq m$ the following property holds:

if $Q_i=\forall$ then $\overline C'_i\sqe C_i$ and if
$Q_i=\exists$ then $C_i\sqe \inte (C'_i)$. Here $\overline B$ is
the closure of $B$ and $\inte (B)$ is the interior of $B$.

If $c\ll c'$ and $W\in\W$ then the formula $\f[c'][W]$ is called a
strong approximation of $\f[c]$. The following lemma is obvious.

\begin{Lm} \label{l1_pos_bound}
Let $\f(\vector x1n)$ be a positive formula of $L_{\theta}$ of the
form (18), $c_1\ll c_2$, be $\f$-regular $m$-tuples of subsets of
$A$, let $W_2\subseteq W_1$ be elements of the uniformity $\W$ and
$\la\A_f,j_f\ra$ be a $(C,W)$-approximation of $\A$ for some $\la
C,W \ra\in M$. Then
\begin{enumerate}
\item $\forall\vector a1n\in A \left(\A\models\f[c_1](\vector
a1n)\Longrightarrow\A\models\f[c_2](\vector a1n)\right)$;

\item $\forall \vector \a1n\in A_f
(\A_f\models\f[c_1][W_2](\vector
\a1n)\Longrightarrow\A_f\models\f[c_1][W_1](\vector \a1n))$;

\item $\forall \vector \a1n\in A_f
(\A_f\models\f[c_1][W_2](\vector
\a1n)\Longrightarrow\A_f\models\f[c_2][W_2](\vector \a1n))$.
\end{enumerate}
\end{Lm}

This notion of approximation for positive bounded formulas is
similar to the one introduced in \cite{he} \cite{hm} \cite{hi})
for structures based on Banach spaces.

The following theorem is the main result of this section.

\begin{Th} \label{pos_bound}
Let $\f[c](\vector x1n)$ be a positive bounded formula and
$\vector a1n\in A$. Then $\A\models \f[c](\vector a1n)$ iff  for
any strong approximation $\f[c'][W']$ of $\f[c]$ there exists a
pair $\langle  C_0,W_0\rangle\in M$ such that the following
conditions hold:

1) $\bigcup\limits_{i=1}^nW_0(a_i)\sqe C_0$;

2) for any $\langle  C,W\rangle\leq \langle  C_0,W_0\rangle$, for
any $(C,W)$-approximation $\langle \A_f,j_f\rangle$ of $\A$ and
for any $\vector b1n\in A_f$ such that $\la a_i,j(b_i)\ra\in W_0,\
i=1,\dots,n$, one has $\A_f\models\f[c'][W'](\vector b1n)$.
\end{Th}

If for some property P there exists a $\la C_0, W_0\ra\in M$ such
that P holds for all $(C,W)$-approximations of $\A$ such that $\la
C,W\ra\leq \la C_0,W_0\ra$, then we say that P holds for all
precise enough approximations of $\A$.

\begin{Cor} \label{pos_bound_sent}
A positive bounded sentence $\f[c]$ holds in $\A$ iff all of its
strong approximations $\f[c'][W]$ hold in all precise enough
approximations of $\A$.
\end{Cor}

>From the point of view of numerical systems implemented in
computers this corollary means that approximate versions of
positive bounded theorems about the reals hold for numerical
computer systems that simulate the field of reals in powerful
enough computers.

Before we start to prove Theorem \ref{pos_bound}, consider the
following three examples. In these examples we deal with the
algebra $\la\R;1,+,\times\ra$ and its $(a,\e)$-approximations
$\la\A_f, j_f\ra$ (see Example 1) such that $j_f$ is the inclusion
map. According to Definition \ref{c-w-appr}(3), in this case we
say that $\A_f$ is an $(a,\e)$-approximation of $\R$.

\bigskip\noindent
{\bf Example 6} Fix any positive $d>1$. Then the following
positive bounded formula holds for the field $\R$:
$$
\forall_D\,x\,\exists_{\overline D}\, y(xy=1),
$$
where $D=\{x\in\R\ |\ d^{-1}< |x|< d\}$.

It is easy to see that for any strong approximation of this
formula there exists a finer strong approximation of the following
form:
$$
\forall_C\, x\,\exists_B\,y(|xy-1|<\d), \eqno (20)
$$
where $C=\{x\in\R\ |\ c^{-1}<|x|<c\},\ B=\{x\in\R\ |\ b^{-1}\leq
|x|\leq b\}$, $1<c<d<b$ and $\d>0$.

We have to show that there exist $a_0,\e_0$ such that for any
$a>a_0,\ \e<\e_0$, formula (20) holds for any finite $(a,\e)$-approximation $\A_f$ of $\R$. Fix any $x$ such that
$c^{-1}<|x|<c$ and let $y=x^{-1},\ b^{-1}<|y|<b$. Take
$\xi,\eta\in A_f$ such that $|x-\xi|<\e$ and $|y-\eta| <\e$. The
$a$ and $\e$ have to satisfy the following conditions:
$\xi,\eta,\xi\times \eta\in [-a,a],\ |\xi\otimes\eta-1|<\d$, where
$\otimes$ is the multiplication in $\A_f$. By the definition of
$(a,\e)$-approximation, it is easy to see that the following
$a_0$ and $\e_0$ satisfy the required conditions:
$$
\e_0=\sqrt{\left(\frac{2b+1}2\right)^2+\d}-\frac{2b+1}2,\
a_0=\max\{b+\e_0,(2b+\e_0)\e_0+1\}.
$$

\bigskip\noindent
{\bf Example 7} Let $\f(x,y)$ be the positive formula $\exists
z(x+z^2=y)$, which defines the relation $\leq$ in $\R$.

Consider a bounded version of this formula
$\f[b](x,y)=\exists_{|z|\leq b}(x+z^2=y)$, which defines the
relation $x\leq y\leq x+b^2$. A strong approximation of this
formula is of the form $\f[c][\a](x,y)=\exists_{|z|\leq
c}(|x+z^2-y|<\a)$ for some $\a>0$ and $0<b<c$. Let $x_0,y_0\in
[-d,d]$ and $\R\models\f[b](x_0,y_0)$. Put $a_0=(c+\a)^2+d+\a+1$
and $\e_0=\max\{c-b,\frac{\a}{5+2a_0}\}$. Then it is easy to see
that $\A_f\models\f[c][\a](\xi,\eta)$, where $\A_f$ is a $(a,\e)$-approximation  of $\R$, $a>a_0,\ \e<\e_0$, $\xi,\eta\in
A_f$ and $|x_0-\xi|<\e_0,\ |y_0-\eta|<\e_0$. If $x_0>y_0$ and
$\a<\frac 12(x_0-y_0)$, we may take $\e_0$ such that
$x_0-\e_0>y_0+\a$. Thus, the formula $\f[c][\a](\xi,\eta)$ fails
in $\A_f$ for any $(a,\e)$-approximation $\A_f$ of $\R$ and
for any $\xi,\eta\in A_f$ such that $|x_0-\xi|<\e_0,\
|y_0-\eta|<\e_0$. A similar consideration holds for $y_0>x_0+b^2$.

\bigskip\noindent
{\bf Example 8} The relation $<$ can also be defined by a positive
formula. Indeed:
$$
x<y\Longleftrightarrow\exists z((y-x)z^2=1)=\f(x,y)
$$
A bounded version of this formula $\f[b](x,y)=\exists_{|z|\leq
b}((y-x)z^2=1)$ defines the relation $y>x+\frac 1{b^2}$. A strong
approximation of this formula is of the form
$\f[c][\a](x,y)=\exists_{|z|\leq c}(|(y-x)z^2-1|<\a)$ for some
$\a>0$ and $0<b<c$. It is easy to see that for $\a<1$, for small
enough $\e$, big enough $a$ and for any $(a,\e)$-approximation $\A_f$ of $\R$ if $\xi,\eta\in [-a,a]$ then
$\A_f\models\f[c][\a](\xi,\eta)\Longleftrightarrow
\eta>\xi+\frac{1-\a}{c^2}$.

\begin{rem}
  It is easy to see that the relation $u \prec v$ between normalized
  floating-point numbers $u,v$, introduced in \cite[page 200]{kn},
  is a special case of the approximation $\f[c][\alpha]$ in Example
  8.
\end{rem}

\begin{rem} \label{quant_free}
By a classical result of Tarski \cite{Tar} any formula in the
signature $\la +,\times\ra$ is equivalent in the first order
theory of the ordered field of real numbers ($\mbox{Th}(\R)$) to a
quantifier free formula in the signature $\la 1,+,\times,\leq\ra$.
Therefore, the examples considered above show that any formula of
the language of rings is equivalent in $\mbox{Th}(\R)$ to a
positive formula and thus has its approximate versions.
\end{rem}

Let the topological space $A$ be totally disconnected; i.e., the
clopen sets form a base of its topology. Consider a positive
bounded formula $\f[c]$ with an $m$-tuple $c$ that consists of
clopen sets. In this case we say that $c$ is clopen. Since for a
clopen set $V$ one has $\overline V\sqe V$ and $V\sqe\inte(V)$,
then for a clopen $m$-tuple $c$ one has $c\ll c$. Thus if
$W\in\W$, then $\f[c][W]$ is a strong approximation of $\f[c]$. So
the formulation of Theorem \ref{pos_bound} can be simplified for
this case.

\begin{Cor} \label{tot_disc}
Let $\A$ be a totally disconnected algebra, $\f[c](\vector x1n)$
be a positive bounded formula (19) with a clopen $m$-tuple $c$,
and $\vector a1n\in A$. Then $\A\models \f[c](\vector a1n)$ iff
for any $W'\in\W$ there exists a pair $\langle  C_0,W_0\rangle\in
M$ such that the following conditions hold:

1) $\bigcup\limits_{i=1}^nW_0(a_i)\sqe C_0$;

2) for any $\langle  C,W\rangle\leq \langle  C_0,W_0\rangle$, for
any $(C,W)$-approximation $\langle \A_f,j_f\rangle$ of $\A$, and
for any $\vector b1n\in A_f$ such that $\la a_i,j(b_i)\ra\in W_0,\
i=1,\dots,n$, one has $\A'\models\f[c][W'](\vector b1n)$.
\end{Cor}

Now we turn to the proof of Theorem \ref{pos_bound}. First we
consider an equivalent nonstandard statement.

Let $\la\A_h,j\ra$ be a hyperfinite approximation of $\A$ in the
sense of Definition \ref{nonst-appr}. Then a strong approximation
$\f[c][W]$ of a positive formula $\f$ in the form (18) has an
obvious interpretation in $\A_h$: a quantifier $Q_Cx\dots$ is
interpreted as on page \pageref{ibq} and a formula $\la
t_1,t_2\ra\in W$ is interpreted by $\la j(t_1),j(t_2)\ra\in\* W$.
Obviously, the statements (2) and (3) of Lemma \ref{l1_pos_bound}
hold for hyperfinite approximations of $\A$.

\begin{Lm} \label{l2_pos_bound}
For any $\vector {\b}1n\in (A_h)_b$ one has \newline
\centerline{$\A_h\models\f[c]_{\sim}(\vector {\b}1n)
\Longleftrightarrow \forall W\in\W_{\l}\
\A_h\models\f[c][W](\vector {\b}1n).$}
\end{Lm}

\medskip\noindent
{\bf Proof} Obviously, $\A_h\models\f[c]_{\sim}(\vector
{\b}1n)\Longrightarrow\forall W\in\W_{\l}\
\A_h\models\f[c][W](\vector {\b}1n).$ So we have to prove only the
converse implication. Consider first the case of a quantifier free
formula, i.e., the case when $\f=\psi$ in the form (18). We have
$\psi=P_1\lor\dots\lor P_r$, where each $P_i$ is a conjunction of
atomic formulas. Assume that $\forall W\in\W_{\l}$ one has
$\psi[W]$. If $\psi_{\sim}$ is false then for each $i\leq r$ there
exists $W_i\in\W_{\l}$ such that $P_i[W_i]$ is false. Take $W\in
W_{\l}$ such that $W\sqe\bigcap\limits_{i=1}^rW_i$. Then by Lemma
\ref{l1_pos_bound}(2) for any $i\leq r$ the formula $P_i[W]$ is
false. Thus, the formula $\psi[W]$ is false.

We have to prove now that
$$
\forall W\in\W_{\l}Q_{1C_1}y_1\dots
Q_{mC_m}y_m\psi[W]\Longrightarrow Q_{1C_1}y_1\dots
Q_{mC_m}y_m\forall W\in\W_{\l}\psi[W]
$$

To prove this implication, it is enough to prove that for any
positive bounded formula $\tau(x)$ and any compact set $C$ one has
$$
\forall W\in\W_{\l}\exists_C x\tau[W](x)\Longrightarrow \exists_C
x\forall W\in\W_{\l}\tau[W](x). \eqno(21)
$$

Assume that the left hand side of this implication holds. Put
$B(W)=\{x\ |\ j(x)\in\* C, \tau[W](x)\}$. Then
$B(W)\neq\emptyset$. Since for any $\vector W1s\in\W_{\l}$ there
exists $W\in\W_{\l}$ such that $W\sqe\bigcap\limits_{i=1}^sW_i$,
using Lemma \ref{l1_pos_bound}(2), we obtain that the family
$\{B(W)\ |\ W\in\W_{\l}\}$ has the finite intersection property.
Thus, by saturation, we obtain that the right hand side of the
implication (21) holds.  \hfill $\Box$

\begin{Lm} \label{l3_pos_bound}
Let $\f(\vector x1n)$ be a positive formula in $L_{\theta}$ of the
form (18), $c=\la\vector C1m\ra$ a $\f$-regular $m$-tuple of
subsets of $A$, $\vector a1n\in A$, and $\la\A_h,j\ra$ a
hyperfinite approximation of $\A$. Then $\A\models\f[c](\vector
a1n)$ iff for any $\f$-regular $c'=\la\vector{C'}1m\ra$ such that
$c'\gg c$ and for any $\vector {\a}1n\in A_h$ such that
$j(\a_i)\approx a_i,\ i=1,\dots,n$ one has
$\A_h\models\f[c']_{\sim}(\vector {\a}1n)$.
\end{Lm}

\medskip\noindent
{\bf Proof} We prove this lemma by induction on $m$. For $m=0$ it
follows from Proposition \ref{nonst_form}. Assume that it is
proved for $m-1$. Denote by $\tau(y_1,\vector x1n)$ the formula
$Q_2y_2\dots Q_my_m\psi(\vector x1n,\vector y1m)$, by $c_-$ the
$(m-1)$-tuple $\la\vector C2m\ra$, by  $c'_-$ the $(m-1)$-tuple
$\la\vector {C'}2m\ra$, so, that $\f=Q_1y_1\tau,\
\f[c]=Q_{1C_1}y_1\tau[c_-],\ \f[c']=Q_{1C'_1}y_1\tau[c'_-]$.
Consider two cases.

a). $Q_1=\exists$. In this case $C_1, C'_1$ are compact sets and
$C_1\sqe\inte{C'_1}$.

$\Rightarrow$ Let $\A\models\f[c](\vector a1n)$. Then there exists
$b\in C_1$ such that $\A\models\tau[c_-](b,\vector a1n)$. Let
$\b\in A_h$ be such that $j(\b)\approx b$. Then, by the induction
assumption, $\A_h\models\tau[c'_-]_{\sim}(\b,\vector{\a}1n)$.
Since $b\in C_1\sqe\*\inte(C'_1)\sqe\*C'_1$, $j(\b)\approx b$ and
$\inte(C'_1)$ is an open set, we have $j(\b)\in\*\inte(C'_1)$.
This proves that $\A_h\models\f[c']_{\sim}(\vector {\a}1n)$.

$\Leftarrow$. Obviously $C_1=\bigcap\{\overline{W(C_1)}\ |\
W\in\W_{\l}\}$. Fix any $V,W\in\W_{\l}$ . By the assumption of the
lemma, for any $\tau$-regular $(m-1)$-tuple $c'_-\gg c_-$ one has
$\A_h\models \ex_{\overline{W(C_1)}}y_1\tau[c'_-]_{\sim}(\vector
{\a}1n)$. Thus, by Lemma \ref{l2_pos_bound}, we have
$B(W,V,c'_-)=\{\b\in j^{-1}(\*\overline{W(C_1)})\ |\
\A_h\models\tau[c'_-][V](\b,\vector {\a}1n)\}\neq\emptyset$. Let
$\Xi$ be the set of all $\tau$-regular $(m-1)$-tuples $c'_-$. Then
it is easy to see that there exists a cofinal subset $\Xi_{\l}$
(i.e., $\forall d\in\Xi\exists d_1\in\Xi_{\l}(d_1\ll d)$) of
cardinality $\l$. By Lemma \ref{l1_pos_bound}(3), if for some $b,
\vector {\a}1n\in (A_h)_b$ for all $d\in\Xi_{\l}$ one has
$\A_h\models\tau[d][V]$, then the same holds for all $d\in\Xi$. It
is easy to see also that, similar to $\W_{\l}$, the family
$\Xi_{\l}$ has the following property: for any
$c^{(1)},\dots,c^{(s)}\in\Xi_{\l}$ there exists a $c'_-\in
\Xi_{\l}$ such that $c'_-\ll c^{(1)},\dots,c'_-\ll c^{(s)}$. All
this shows that the family $\{B(W,V,c'_-)\ |\
V,W\in\W_{\l},c'_-\in\Xi_{\l}\}$ has the finite intersection
property and thus, by saturation, has nonempty intersection. By
our construction and Lemma \ref{l2_pos_bound}, any element $\b$ in
this intersection has the following properties: $j(\b)\in\* C_1$
and $\A_h\models\tau[c'_-]_{\sim}(\b,\vector {\a}1n)$ for any
$c'_-\in\Xi$. By the induction assumption this implies that
$\A\models\tau[c_-](\o j(\b),\vector a1n)$. Since $C$ is a compact
set, we obtain that $\o j(\b)\in C$. This proves a).

b)$Q_1=\forall$. In this case $C_1$ and $C'_1$ are relatively
compact open sets and $\overline{C'_1}\sqe C_1$

$\Rightarrow$ Let $\A\models\f[c](\vector a1n)$ and $\vector
{\a}1n\in (A_h)_b$ be such that $j(\a_i)\approx a_i$. Take any
$\b\in A_h$ such that $j(\b)\in\* C'_1$. Then $b=\o j(\b)\in C_1$
and, thus, $\A\models\tau[c_-](b,\vector a1n)$. By the induction
assumption, $\A_h\models\tau[c'_-]_{\sim}(\b,\vector {\a}1n)$.
This completes the proof.

$\Leftarrow$.  Let $b\in C_1$. Obviously, there exists an open $D$
such that $\overline D\sqe C_1$ and $b\in D$. Since for any
$c'_-\gg c_-$ and for any $\vector {\a}1n\in (A_h)_b$ such that
$j(\a_i)\approx a_i$ one has
$\A_h\models\forall_Dy_1\tau[c'_-]_{\sim}(y_1,\vector {\a}1n)$, we
obtain that for any $\b\in A_h$ such that $j(\b)\approx b$ one has
$\A_h\models\tau[c']_{\sim}(\b,\vector {\a}1n)$. Thus, by the
induction assumption $\A\models\tau[c](b,\vector a1n)$. \hfill
$\Box$

\medskip\noindent
{\bf Proof of Theorem \ref{pos_bound}} ($\Rightarrow$) Let
$\A\models\f[c](\vector a1n)$. Fix a strong approximation
$\f[c'][W'](\vector x1n)$ of $\f[c](\vector x1n)$ . Consider the
internal set $N$ of all pairs $\la C_0,W_0\ra\in\* M_{\l}$ such
that for all $\la C,W\ra\leq\la C_0,W_0\ra$, for any $(
\*C,\*W)$-approximation $\la\A_h,j\ra$ of $\*\A$ and for any
$\vector {\a}1n\in A_h$ satisfying the condition $\la
a_i,j(\a_i)\ra\in W_0$ one has $\A_h\models\f[c][W](\vector
{\a}1n)$. Lemmas \ref{l2_pos_bound}, \ref{l3_pos_bound} and
\ref{infin-appr} imply that $N$ contains all infinitesimal pairs
$\la C_0,W_0\ra\in\* M_{\l}$. By Lemma \ref{about_M}(2), there
exists $\la C_0,W_0\ra\in M_{\l}$ such that $\la\* C_0,\*
W_0\ra\in N$. By the transfer principle, this completes the proof.

$\Leftarrow$. Let $c'\gg c$, $\la\A_h,j\ra$ be a hyperfinite
approximation of $\A$ and $\vector {\a}1n\in A_h$ be such that
$j(\a_i)\approx a_i$. Then for any $\la C,W\ra\in M_{\l}$ the pair
$\la A_h,j\ra$ is a $(\* C,\* W)$-approximation of $\*\A$ and $\la
j(\a_i),a_i\ra\in\* W,\ i=1,\dots,n$. Thus, by the conditions of
the theorem, $\A_h\models\f[c'][W'](\vector {\a}1n)$ for any
$W'\in\W_{\l}$. By Lemma \ref{l2_pos_bound},
$\A_h\models\f[c]_{\sim}(\vector {\a}1n)$. Hence, by Lemma
\ref{l3_pos_bound} one has $\A\models\f[c](\vector a1n)$. \hfill
$\Box$

The following corollary of Theorem \ref{pos_bound} shows that the
approximation of continuous functions by polynomials on closed
intervals holds for all precise enough approximations of the field
$\R$ (cf. the example concerning $\sin x$ which was discussed in
the Introduction).

\begin{Cor} \label{poly}
Let $\A=\la\R,\s\ra$ be such that $\s$ contains the symbols $+$ and
$\times$ and a unary function symbol $g$. Suppose that the
continuous function $g$ is approximable on an interval $[-d,d]$ by a
polynomial $b_nx^n+\dots +b_1x+b_0$ with accuracy $\d$. Then for
any $0<d'<d$ and $\d'>\d$ there exist $a_0,\e_0>0$ such that
any $(a,\e)$-approximation
$\la\A_f,j_f\ra$ of $\A$ with $a>a_0$ and $\e<\e_0$ has the following property:
for any $\vector {\b}0n, \xi \in A_f$, if
$|j(\b_i)-b_i|<\e_0$ for all $i=0,\dots, n$ and
$j(\xi)\in [-d',d']$, then
$|j_f(g_f(\xi))-j_f(\b_n\xi^n+\dots +\b_1\xi+\b_0)|<\d'$.
(Here $g_f$ is the interpretation of the symbol $g$ in $\A_f$.)
\end{Cor}


\bigskip\bigskip

Universidad Autonoma de San Luis Potosi, Mexico (Glebsky)

Eastern Illinois University, USA (Gordon)

University of Illinois at Urbana-Champaign, USA (Henson)

\bigskip\bigskip

1991 {\it Mathematics Subject Classification}. Primary 26E35,
03H05; Secondary 28E05, 42A38

\newpage

Instituto de Investigacion en Communicacion Optica

Universidad Autonoma de San Luis Potosi, Mexico

AvKarakorum 1470

Lomas 4ta Session

San Luis Potosi SLP 7820

Mexico

e-mail:glebsky@cactus.iico.uaslp.mx

\bigskip

Department of Mathematics and Computer Science

Eastern Illinois University

600 Lincoln Avenue

Charleston, IL 61920-3099

USA

e-mail: cfyig@eiu.edu

\bigskip

Department of Mathematics

University of Illinois at Urbana-Champaign

1409 West Green Street

Urbana, IL 61801

USA

www: http://www.math.uiuc.edu/$\sim$henson

\end{document}